\documentclass[10pt]{article}
\usepackage[a4paper, total={7in, 8in}]{geometry}
\usepackage{url,floatflt}
\usepackage{helvet,times}
\usepackage{graphics}
\usepackage{mathptmx,amsmath,amssymb,bm}
\usepackage{float}
\usepackage[bf,hypcap]{caption}
\usepackage{enumerate}

\title{An initial-boundary problem for a mixed fractional wave equation}
\author{Erkinjon Karimov{$^{2,3}$}, Nasser Al-Salti{$^{1, 3,*}$}, Muna Al-Ghabsi{$^{1}$}\\
1: Department of Applied Mathematics and Science,\\ National University of Science and Technology, Muscat, Oman \\ 2: Ghent University, Ghent-9000, Belgium \\ 3: FracDiff Research Group (DR/RG/03), 
Sultan Qaboos University, Muscat, Oman
}

\begin{document}

\maketitle{}

\begin{abstract}
We aim to prove a unique solvability of an initial-boundary value problem (IBVP) for a time-fractional wave equation in a rectangular domain. We exploit the spectral expansion method as the main tool and used the solution to Cauchy problems for fractional-order differential equations. Moreover, we apply certain properties of the Mittag-Leffler-type functions of single and two variables to prove the uniform convergence of the solution to the considered problem, represented in the form of infinite series.
\end{abstract}

{\bf Keywords:\,}{Wave equation, Fractional Differential Equations, Mittag-Leffler type functions, Cauchy problem, mixed type equation.


\section{\; Introduction}

It is well known that partial differential equations (PDEs) are fundamental in mathematical modeling as they describe how physical, biological, and engineering systems evolve over time and space (see \cite{R1} for its general theory and \cite{R2}, \cite{R3}, \cite{R4} for applications). We will direct our attention to hyperbolic type equations \cite{R5}, which are linked to many topics such as waves, compressible fluid flow, supersonic flows, etc. \cite{R6}.

In \cite{R16}, author considered a fractional wave equation that contains fractional derivatives of the same order $\alpha$ ($1\le\alpha\le 2$), both in space and in time. Different velocities of damped waves were shown to be described by the fundamental solution of the fractional wave equation, and they (velocities) will be constant depending on the order of the equation $\alpha$. Another application of fractional wave equations is related to the well-known telegraph process \cite{R17}. We note recent publications \cite{R18}-\cite{R20}, in which fractional telegraph equations were targeted at different view points. Detailed information on fractional wave equations can be found in \cite{R21}.

Acoustic problems and gas flow in tubes can be described via boundary-value problems for hyperbolic wave equations. If a damping force is acting or some barriers are involved in acoustics for some time and then disappear, the mathematical model of such a process should be described via mixed-wave equations. During the damped phase, the wave loses energy due to damping, resulting in a decrease in amplitude, and after the damping is removed, the wave propagates freely without further energy loss, maintaining its amplitude. This motivated us to formulate a well-posed problem for such a situation and investigate its unique solvability. Hence, this paper considers an initial boundary problem for a mixed fractional wave equation. We note that the combination of two wave equations was considered in recent work \cite{R15}. Namely, aime-fractional telegraph and one classical telegraph equation were considered in the pentagonal dodomain, imposing non-local conditions. 

The remainder of the paper is organized as follows. First, we present some preliminaries on fractional derivatives, Mittag-Leffler type functions, and their properties and solutions to some related Cauchy problems. Then, we present our main work on problem formulation and formal solution. 

\section{\; Preliminaries}
This section presents definitions of fractional derivatives, Mittag-Leffler type functions, and their properties. We also present solutions to two Cauchy problems that will be used during this study.
\subsection{\; Riemann-Liouville and Caputo Fractional Operators}
{\bf Definition 1.} \cite{R7} 
	For $\alpha \notin \mathbb{N}\cup \left\{ 0 \right\}$, the Riemann-Liouville fractional integral $I_{ax}^\alpha f$ of order $\alpha\in \mathbb{C}$ is defined by
	\begin{equation}\label{RL_int}
		I_{ax}^\alpha f(x)=\frac{1}{\Gamma(\alpha)}\int\limits_a^x \frac{f(t) \,dt}{(x-t)^{1-\alpha}}, \,\, \quad x>a,\,\Re(\alpha)>0
	\end{equation}
	and the Riemann-Liouville and the Caputo fractional derivatives of order $\alpha$ are defined by: 
	\begin{equation}\label{RL_der}
		^{RL}D_{ax}^\alpha y(x)=\frac{d^n}{dx^n}I_{ax}^{n-\alpha} y(x)=\frac{1}{\Gamma(n-\alpha)}\frac{d^n}{dx^n}\int\limits_a^x \frac{y(t) \, dt}{(x-t)^{\alpha-n+1}}, \,\quad x>a,\,n=[\Re(\alpha)]+1,
	\end{equation}
	\begin{equation}\label{C_der}
		^{C}D_{ax}^{\alpha}y(x)=\frac{1}{\Gamma \left( n-\alpha  \right)}\int\limits_{a}^{x}{\frac{{{y}^{\left( n \right)}}\left( t \right) \, dt}{{{\left( x-t \right)}^{\alpha -n+1}}}}, \,\, \quad x>a,\, n=[\Re(\alpha)]+1,
	\end{equation}
respectively.

\subsection{\; Mittag-Leffler type functions}

The two-parameter Mittag-Leffler function is defined by a series expansion as follows:
{\bf Definition 2.} \cite{R8} For $\Re(\alpha)>0, \beta \in \mathbb{C}$, the two-parameter Mittag-Leffler function is defined as
	\begin{equation}\label{ML_2}
		{{E}_{\alpha ,\beta }}\left( z \right)=\sum\limits_{k=0}^{\infty }{\frac{{{z}^{k}}}{\Gamma \left( \alpha k+\beta  \right)}}
	\end{equation}
and has the following properties:
\begin{enumerate}[(i)]
	\item 	\quad ${{E}_{\alpha ,\beta }}\left( z \right)=\frac{1}{\Gamma \left( \beta  \right)} + z{{E}_{\alpha ,\alpha +\beta }}\left( z \right)$;
	\medskip
	\item \quad $\left(\dfrac{d}{dz}\right)^m \left[z^{\beta -1} E_{\alpha,\beta}\left(z^\alpha \right)\right]={{z}^{\beta- m -1}}E_{\alpha ,\beta -m }\left({{z}^{\alpha }} \right) \quad m \in \mathbb{N}$;
	\medskip
	\item \quad $^{RL} D_{0t}^\alpha \left(t^{\gamma -1} E_{\beta, \gamma}\left(at^\beta \right)\right)=t^{\gamma-\alpha-1} E_{\beta,\gamma-\alpha}\left(at^\beta \right), \quad \Re(\alpha)>0, \Re(\beta)>0, \Re(\gamma)>0, a \in \mathbb{R}$ ;
	\medskip
	\item \quad $I_{at}^\alpha \left(t^{\gamma -1} E_{\beta, \gamma}\left(at^\beta \right)\right)=t^{\alpha+\gamma-1} E_{\beta,\alpha+\gamma}\left(at^\beta \right), \quad \Re(\alpha)>0, \Re(\beta)>0, \Re(\gamma)>0, a \in \mathbb{R}$ ;
\end{enumerate}

{\bf Lemma 1.} \cite{R9} If $0<\alpha<2, \beta\in \mathbb{R}, \dfrac{\pi \alpha}{2}< \mu < min \{\pi, \pi \alpha \}$, then
$$ | E_{\alpha, \beta}(z)|\le \dfrac{c}{1+|z|}, \quad \mu\le|arg(z)|\le \pi, |z|\ge 0.$$

Now, we present the multi-variable Mittag-Leffler function in a particular case as follows:

{\bf Definition 3.} \cite{R10} For $\rho>0$ and $\alpha,\,\beta\in\mathbb{R}^+$, $|x|<\infty,\,|y|<\infty$
	\begin{equation}\label{ML_mv}
		E_{(\alpha, \beta), \rho} (x,y)=\sum\limits_{n=0}^{\infty} \, \sum\limits_{k=0}^{n} \dfrac{n!}{k!(n-k)!} \, \dfrac{x^k y^{n-k}}{\Gamma(\rho+\alpha n + \beta k)}
	\end{equation} 
	presents the bi-variate case of multi-variable Mittag-Leffler function.

We note the following properties of this function:

\begin{enumerate} [(I)]
	\item \quad 
	$$
	\begin{array}{l}
		m_1 t^{\alpha-\beta} E_{(\alpha-\beta, \alpha), \alpha-\beta+\rho} \left(m_1t^{\alpha -\beta}, m_2 t^{\alpha}\right) +m_2 t^{\alpha} E_{(\alpha-\beta, \alpha), \alpha+\rho} \left(m_1t^{\alpha -\beta}, m_2 t^{\alpha}\right)=\\
		E_{(\alpha-\beta, \alpha), \rho} \left(m_1t^{\alpha -\beta}, m_2 t^{\alpha}\right)- \dfrac{1}{\Gamma(\rho)};
	\end{array}
	$$
	\medskip
	\item \quad $ \dfrac{d}{dt} \left(t^{\alpha} E_{(\alpha-\beta, \alpha), \alpha+1} \left(m_1t^{\alpha -\beta}, m_2 t^{\alpha}\right) \right)=t^{\alpha-1} E_{(\alpha-\beta, \alpha), \alpha} \left(m_1t^{\alpha -\beta}, m_2 t^{\alpha}\right)$;
	\medskip
	\item {\quad $ \int \limits_{0}^{t} \left(z^{\alpha-1} E_{(\alpha-\beta, \alpha), \alpha} \left(m_1z^{\alpha -\beta}, m_2 z^{\alpha}\right) \right)=t^{\alpha} E_{(\alpha-\beta, \alpha), \alpha+1} \left(m_1t^{\alpha -\beta}, m_2 t^{\alpha}\right)$.}
\end{enumerate}

We will use the following statement, which follows from Lemma 5 of \cite{R11}:

{\bf Lemma 2.}
Let $\rho>0$ and $\alpha>\beta>0$ be given with $\alpha\in(0,2)$. Assume that $x<0$ and there exist $K>0$ such that $-K\le y<0$. Then there exists a constant $C>0$ depending only on $K,\,\alpha,\,\beta$ and $\rho$ such that
\[
\left|E_{(\alpha,\beta),\,\rho}(x,y)\right|\le \frac{C}{1+|x|}.
\]

\subsection{\; Cauchy Problems}
Here we present two Cauchy problems. First, we consider the following Cauchy problem:
\begin{equation}\label{CP1}
	\left\{	\begin{array}{ll}
		^CD_{0t}^{\alpha_1} y(t)+ \mu ^CD_{0t}^{\alpha_2} y(t)+\lambda y(t)=f(t), \quad t>0 \vspace{0.2cm}\\
		y(0)=A, \quad y'(0)=B, \quad 0<\alpha_2< 1, \quad 1<\alpha_1< 2, \quad \lambda, \mu \in \mathbb{R},
	\end{array} \right.
\end{equation}
which has the following solution that is expressed in terms of the multi-variable Mittag-Leffler function \cite{R10}:
\begin{eqnarray} \label{CP1_sol}
	y(t)&=& A \left[1-\mu t^{\alpha_1 - \alpha_2} E_{(\alpha_1-\alpha_2, \alpha_1),\alpha_1-\alpha_2+1}\left(-\mu t^{\alpha_1-\alpha_2}, -\lambda t^{\alpha_1} \right)\right] \nonumber\\ 
	&+& B \left[1-\mu t^{\alpha_1 - \alpha_2} E_{(\alpha_1-\alpha_2, \alpha_1),\alpha_1-\alpha_2+2}\left(-\mu t^{\alpha_1-\alpha_2}, -\lambda t^{\alpha_1} \right) -\lambda t^{\alpha_1} E_{(\alpha_1-\alpha_2, \alpha_1),\alpha_1+2}\left(-\mu t^{\alpha_1-\alpha_2}, -\lambda t^{\alpha_1} \right)\right]\\
	&+& \int \limits_{0}^{t} s^{\alpha_1 - 1} E_{(\alpha_1-\alpha_2, \alpha_1),\alpha_1}\left(-\mu s^{\alpha_1-\alpha_2}, -\lambda s^{\alpha_1} \right) f(s) \, ds \nonumber
\end{eqnarray}
Now, we consider the following second Cauchy problem:
\begin{equation} \label{CP2}
	\left\{		\begin{array}{ll}
		^CD_{at}^{\beta} y(t)+\lambda y(t)=f(t), \quad t>a  \vspace{0.2cm}\\
		{	y(a)=C, \quad y'(a)=G,} \quad 1<\beta< 2, \quad \lambda, \mu \in \mathbb{R},
	\end{array} \right. 
\end{equation}
along with its solution, which is given by \cite{R7}
\begin{equation} \label{CP2_sol}
	y(t)= CE_{\beta,1}\left(-\lambda (t-a)^{\beta}\right)+ G(t-a) E_{\beta,2}\left(-\lambda (t-a)^{\beta}\right) + \int \limits_{a}^{t} f(s)(s-a)^{\beta - 1} E_{\beta, \beta}\left(-\lambda (s-a)^\beta \right) \, ds. 
\end{equation}
We are now ready to formulate our main problem and investigate its solvability. This will be done in the next section.
\vskip 2mm

\section{\; Problem Formulation and Formal Solution}

Our aim is to find a function $u(t,x) \in C(\bar{\Omega}) \cap C_x^2(\Omega), u(.,x) \in AC^2[0,a]\cup[a,b],$ which satisfies the equation:
\begin{equation} \label{u_eqn}
	u_{xx}(t,x)+f(t,x)=\left\{  
	\begin{array}{ll}
		^CD_{0t}^{\alpha_1} u(t,x) + \mu ^CD_{0t}^{\alpha_2} u(t,x), \quad (t, x) \in \Omega_1 \vspace{0.2cm}\\
		^CD_{at}^\beta u(t,x), \hspace{2.65cm} (t,x) \in \Omega_2
	\end{array}
	\right.
\end{equation}
together with the initial condition:
\begin{equation}\label{u_initial}
	u(0,x)=\phi(x), \quad 0\le x\le 1,
\end{equation}
boundary conditions:
\begin{equation}\label{bc}
	u(t,0)=0,\,\,u(t,1)=0,\,\,\,0\le t\le b,   
\end{equation}
the terminal condition:
\begin{equation}\label{u_terminal}
	u(b,x)=\psi(x), \quad 0\le x\le 1,
\end{equation}
and the transmitting conditions:
\begin{equation} \label{u_transm1}
	u(a+0,x)=u(a-0,x), \quad 0\le x\le 1,
\end{equation}
\begin{equation} \label{u_transm2}
	\lim\limits_{t\rightarrow \, a^+} {^C}D_{at}^\beta u(t,x)=u_t(a-0,x), \quad 0\le x\le 1,
\end{equation}
where $\phi(x,), \psi(x), f(t,x)$ are given functions and 
\[ \Omega_1= \left\{(t,x): 0<x<1, \, 0<t<a\right\}, \quad \Omega_2= \left\{(t,x): 0<x<1, \, a<t<b \right\}, \quad \Omega=\Omega_1\cup \Omega_2\]
with $a, b, \mu, \alpha_1, \alpha_2, \beta \in \mathbb{R}$ such that $b>a>0, \, 0<\alpha_2<1, \, 1<\alpha_1, \beta < 2$.

In particular, we look for a solution to the IBVP (\ref{u_eqn})-(\ref{u_transm2}) in the form of a series expansion using an orthogonal basis which can be obtained by considering the corresponding homogeneous problem and using the method of separation of variables. This will lead to the following Sturm-Liouville problem:
\begin{equation} \label{SLP}
	X''(x)+\lambda X(x)=0, \quad X(0)=x(1)=0,
\end{equation}
whose eigenvalues are given by $\lambda_{n}=(n\pi)^2, \, n\in \mathbb{N}$ and the corresponding eigenfunctions are 
\begin{equation}
	X_n(x)=\sin n\pi x,
\end{equation}
which form an orthonormal basis in $L_2(0,1)$. Hence, the solution $u(t,x)$ and the given function $f(t,x)$ can be represented as follow:

\begin{equation} \label{u_FS}
	u(t,x)=\sum \limits_{n=1}^\infty T_n(t) \sin n\pi x, \quad 0\le t\le b,
\end{equation}

\begin{equation} \label{f_FS}
	f(t,x)=\sum \limits_{n=1}^\infty f_n(t) \sin n\pi x, \quad 0\le t\le b,
\end{equation}
where $T_n(t)$ are unknown functions to be found and 
$$
f_n(t)=2 \int_{0}^{1} f(t,x) \sin n\pi x \, dx
$$ 
are the Fourier coefficient of $f(t,x)$. Substituting (\ref{u_FS}) and (\ref{f_FS}) into (\ref{u_eqn}), we obtain the following fractional differential equations for $T_n(t)$:
\begin{equation} \label{T_eqn1}
	^CD_{0t}^{\alpha_1}T_n(t) + \mu ^CD_{0t}^{\alpha_2}T_n(t)+\lambda_{n} T_n(t)=f_n(t), \quad 0\le t< a,
\end{equation}
\begin{equation} \label{T_eqn2}
	^CD_{at}^{\beta}T_n(t)+\lambda_{n} T_n(t)=f_n(t), \hspace{2.3cm} a\le t< b,
\end{equation}
and substituting (\ref{u_FS}) and (\ref{f_FS}) into the initial and terminal conditions (\ref{u_initial}) and (\ref{u_terminal}), we will get the following conditions for $T_n(t)$:
\begin{equation}\label{Tn0}
	T_n(0)=\phi_n,
\end{equation}
\begin{equation}\label{Tnb}
	T_n(b)=\psi_n,
\end{equation}
where 
$$\phi_n(t)=2 \int_{0}^{1} \phi(x) \sin n\pi x \, dx,\,\,\,\psi_n(t)=2 \int_{0}^{1} \psi(x) \sin n\pi x \, dx.
$$
Using the condition (\ref{Tn0}) and based on (\ref{CP1}) and (\ref{CP1_sol}), we obtain the following solution to equation (\ref{T_eqn1}):

\begin{eqnarray} \label{Tn1_sol}
	T_n(t)&=& \phi_n \left[1-\mu t^{\alpha_1 - \alpha_2} E_{(\alpha_1-\alpha_2, \alpha_1),\alpha_1-\alpha_2+1}\left(-\mu t^{\alpha_1-\alpha_2}, -\lambda_n t^{\alpha_1} \right)\right] \nonumber\\ 
	&+& B_n \left[1-\mu t^{\alpha_1 - \alpha_2} E_{(\alpha_1-\alpha_2, \alpha_1),\alpha_1-\alpha_2+2}\left(-\mu t^{\alpha_1-\alpha_2}, -\lambda_n t^{\alpha_1} \right) -\lambda_n t^{\alpha_1} E_{(\alpha_1-\alpha_2, \alpha_1),\alpha_1+2}\left(-\mu t^{\alpha_1-\alpha_2}, -\lambda_n t^{\alpha_1} \right)\right]\\
	&+& \int \limits_{0}^{t} f_n(s) \, s^{\alpha_1 - 1} E_{(\alpha_1-\alpha_2, \alpha_1),\alpha_1}\left(-\mu s^{\alpha_1-\alpha_2}, -\lambda_n s^{\alpha_1} \right)  \, ds, \quad 0\le t \le a, \nonumber
\end{eqnarray}
where $B_n$ {are unknown coefficients to be found}. Similarly, based on (\ref{CP2}) and (\ref{CP2_sol}), we obtain the following solution to equation (\ref{T_eqn2}):
\begin{equation} \label{Tn2_sol}
	T_n(t)= C_nE_{\beta,1}\left(-\lambda_n (t-a)^{\beta}\right)+ G_n(t-a) E_{\beta,2}\left(-\lambda_n (t-a)^{\beta}\right) + \int \limits_{a}^{t} f_n(s)(s-a)^{\beta - 1} E_{\beta, \beta}\left(-\lambda_n (s-a)^\beta \right) \, ds, \,\, a\le t\le b,
\end{equation}
{where $C_n, G_n$ are unknown constants that need to be found.} Hence, in total we have three unknowns, which can be found using the remaining three conditions, namely the two transmitting (\ref{u_transm1}), (\ref{u_transm2}) and the condition (\ref{Tnb}).  Applying condition (\ref{Tnb}) in (\ref{Tn2_sol}), we get 
\begin{equation} \label{UC_eq1}
	\psi_n=C_nE_{\beta,1}\left(-\lambda_n (b-a)^{\beta}\right)+ G_n(b-a) E_{\beta,2}\left(-\lambda_n (b-a)^{\beta}\right) + \int \limits_{a}^{b} f_n(s)(s-a)^{\beta - 1} E_{\beta, \beta}\left(-\lambda_n (s-a)^\beta \right) \, ds.
\end{equation}
Then, using the transmitting conditions (\ref{u_transm1}), which is equivalent to $T_n(a-0)=T_n(a+0)$, we get

\begin{eqnarray} \label{UC_eq2}
	C_n&=& \phi_n \left[1-\mu a^{\alpha_1 - \alpha_2} E_{(\alpha_1-\alpha_2, \alpha_1),\alpha_1-\alpha_2+1}\left(-\mu a^{\alpha_1-\alpha_2}, -\lambda_n a^{\alpha_1} \right)\right] \nonumber\\ 
	&+& B_n \left[1-\mu a^{\alpha_1 - \alpha_2} E_{(\alpha_1-\alpha_2, \alpha_1),\alpha_1-\alpha_2+2}\left(-\mu a^{\alpha_1-\alpha_2}, -\lambda_n a^{\alpha_1} \right) -\lambda_n a^{\alpha_1} E_{(\alpha_1-\alpha_2, \alpha_1),\alpha_1+2}\left(-\mu a^{\alpha_1-\alpha_2}, -\lambda_n a^{\alpha_1} \right)\right]\\
	&+& \int \limits_{0}^{a} f_n(s) \, s^{\alpha_1 - 1} E_{(\alpha_1-\alpha_2, \alpha_1),\alpha_1}\left(-\mu s^{\alpha_1-\alpha_2}, -\lambda_n s^{\alpha_1} \right)  \, ds. \nonumber
\end{eqnarray}
To use the second transmitting condition (\ref{u_transm2}), which leads to \begin{equation}\label{trcond2}	
	\lim\limits_{t\rightarrow \, a^+} {^C}D_{at}^\beta T_n(t)=T'_n(a-0),
\end{equation}
we first need to calculate $T'_n(t)$ for $0<t<a$ and ${^C}D_{at}^\beta T_n(t)$ for $a<t<b$. Hence, differentiating (\ref{Tn1_sol}), we get

\begin{eqnarray} \label{Tn1_der}
	T'_n(t)&=& -\mu \phi_n \dfrac{d}{dt}\left[t^{\alpha_1 - \alpha_2} E_{(\alpha_1-\alpha_2, \alpha_1),\alpha_1-\alpha_2+1}\left(-\mu t^{\alpha_1-\alpha_2}, -\lambda_n t^{\alpha_1} \right)\right] \nonumber\\ 
	&-& \mu B_n \dfrac{d}{dt} \left[ t^{\alpha_1 - \alpha_2} E_{(\alpha_1-\alpha_2, \alpha_1),\alpha_1-\alpha_2+2}\left(-\mu t^{\alpha_1-\alpha_2}, -\lambda_n t^{\alpha_1} \right)\right] -\lambda_n B_n \dfrac{d}{dt} \left[t^{\alpha_1} E_{(\alpha_1-\alpha_2, \alpha_1),\alpha_1+2}\left(-\mu t^{\alpha_1-\alpha_2}, -\lambda_n t^{\alpha_1} \right)\right] \nonumber\\
	&+& \dfrac{d}{dt}\left[\int \limits_{0}^{t} f_n(s) \, s^{\alpha_1 - 1} E_{(\alpha_1-\alpha_2, \alpha_1),\alpha_1}\left(-\mu s^{\alpha_1-\alpha_2}, -\lambda_n s^{\alpha_1} \right)  \, ds, \right], \quad 0<t<a. 
\end{eqnarray}
Using property (II) of the Mittag-Leffler function (\ref{ML_mv}), we get
$$ \dfrac{d}{dt}\left[t^{\alpha_1 - \alpha_2} E_{(\alpha_1-\alpha_2, \alpha_1),\alpha_1-\alpha_2+1}\left(-\mu t^{\alpha_1-\alpha_2}, -\lambda_n t^{\alpha_1} \right)\right]=  t^{\alpha_1 - \alpha_2-1} E_{(\alpha_1-\alpha_2, \alpha_1),\alpha_1-\alpha_2}\left(-\mu t^{\alpha_1-\alpha_2}, -\lambda_n t^{\alpha_1} \right),$$

$$ \dfrac{d}{dt}\left[t^{\alpha_1-\alpha_2} E_{(\alpha_1-\alpha_2, \alpha_1),\alpha_1-\alpha_2+2}\left(-\mu t^{\alpha_1-\alpha_2}, -\lambda_n t^{\alpha_1} \right)\right]=  t^{\alpha_1 -\alpha_2-1} E_{(\alpha_1-\alpha_2, \alpha_1),\alpha_1-\alpha_2+1}\left(-\mu t^{\alpha_1-\alpha_2}, -\lambda_n t^{\alpha_1} \right),$$

$$
\frac{d}{dt}\left[t^{\alpha_1}E_{(\alpha_1-\alpha_2,\alpha_1),\alpha_1+2}\left(-\mu t^{\alpha_1-\alpha_2},-\lambda_nt^{\alpha_1}\right)\right]=t^{\alpha_1-1}E_{(\alpha_1-\alpha_2,\alpha_1),\alpha_1+1}\left(-\mu t^{\alpha_1-\alpha_2},-\lambda_nt^{\alpha_1}\right).
$$
Then using property (I) one can rewrite terms with $B_n$ in other form and considering above-given equalities, (\ref{Tn1_der}) reduces to
\begin{eqnarray}\label{Tn-prime}
	T_n'(t)&=&-\mu\phi_nt^{\alpha_1-\alpha_2-1}E_{(\alpha_1-\alpha_2,\alpha_1),\alpha_1-\alpha_2}\left(-\mu t^{\alpha_1-\alpha_2},-\lambda_nt^{\alpha_1}\right)
	+\frac{B_n}{t}\left[E_{(\alpha_1-\alpha_2,\alpha_1),1}\left(-\mu t^{\alpha_1-\alpha_2},-\lambda_nt^{\alpha_1}\right)-1\right]\nonumber\\
	&+&f_n(t)t^{\alpha_1-1}E_{(\alpha_1-\alpha_2,\alpha_1),\alpha_1}\left(-\mu t^{\alpha_1-\alpha_2},-\lambda_nt^{\alpha_1}\right),\,\,0< t< a.
\end{eqnarray}
Now, applying the operator ${}^CD_{at}^\beta(\cdot)$ to \eqref{Tn2_sol} and we get
\begin{eqnarray*}
	{}^CD_{at}^\beta T_n(t)&=C_n {}^CD_{at}^\beta\left[E_{\beta,1}\left(-\lambda_n(t-a)^\beta\right)\right]+ G_n {}^CD_{at}^\beta\left[(t-a)E_{\beta,2}\left(-\lambda_n(t-a)^\beta\right)\right]\\
	&+{}^CD_{at}^\beta\left[\int\limits_a^tf_n(s)(s-a)^{\beta-1}E_{\beta,\beta}\left(-\lambda_n(s-a)^\beta\right)ds\right].
\end{eqnarray*}
One can evaluate these derivatives using the property (iv) of the function \eqref{ML_2}, but there is a shortcut to get the final result. Precisely, since the equation
$$
{}^CD_{at}^\beta T_n(t)+\lambda_n T_n(t)=f_n(t)
$$
is valid when $t\rightarrow a+0$, then
$$
\lim\limits_{t\rightarrow a+0} {}^CD_{at}^\beta T_n(t)=\lim\limits_{t\rightarrow a+0} \left[f_n(t)-\lambda_n T_n(t)\right].
$$
Hence,
$$
\lim\limits_{t\rightarrow a+0} {}^CD_{at}^\beta T_n(t)=f_n(a)-\lambda_n T_n(a).
$$
If one considers \eqref{Tn2_sol}, it is easy to deduce that
$$
\lim\limits_{t\rightarrow a+0} {}^CD_{at}^\beta T_n(t)=f_n(a)-\lambda_nC_n.
$$
Therefore, \eqref{trcond2} can be rewritten as follows:
\begin{eqnarray}\label{UC_eq3}
	f_n(a)&-&\lambda_nC_n=-\mu\phi_na^{\alpha_1-\alpha_2-1}E_{(\alpha_1-\alpha_2,\alpha_1),\alpha_1-\alpha_2}\left(-\mu a^{\alpha_1-\alpha_2},-\lambda_na^{\alpha_1}\right)\nonumber\\
	&+&\frac{B_n}{a}\left[E_{(\alpha_1-\alpha_2,\alpha_1),1}\left(-\mu a^{\alpha_1-\alpha_2},-\lambda_na^{\alpha_1}\right)-1\right]
	+f_n(a)a^{\alpha_1-1}E_{(\alpha_1-\alpha_2,\alpha_1),\alpha_1}\left(-\mu a^{\alpha_1-\alpha_2},-\lambda_na^{\alpha_1}\right).
\end{eqnarray}
Now, we have \eqref{UC_eq1}, \eqref{UC_eq2} and \eqref{UC_eq3} in order to find the unknown constants $B_n,\,C_n$ and $D_n$. In \eqref{UC_eq1} and \eqref{UC_eq3} instead of $C_n$ we substitute its representation \eqref{UC_eq2}. In that case, from \eqref{UC_eq3} one can find $B_n$ as follows:
$$
B_n=1/\Delta_n\left\{\mu \phi_n a^{\alpha_1-\alpha_2}E_{(\alpha_1-\alpha_2,\alpha_1),\alpha_1-\alpha_2}\left(-\mu a^{\alpha_1-\alpha_2},-\lambda_na^{\alpha_1}\right)+\lambda_n\phi_n [a-\mu a^{\alpha_1-\alpha_2+1}E_{(\alpha_1-\alpha_2,\alpha_1),\alpha_1-\alpha_2+1}\left(-\mu a^{\alpha_1-\alpha_2},-\lambda_na^{\alpha_1}\right)]\right.
$$
$$
+f_n(a)\left[a-a^{\alpha_1+1}E_{(\alpha_1-\alpha_2,\alpha_1),\alpha_1}\left(-\mu a^{\alpha_1-\alpha_2},-\lambda_na^{\alpha_1}\right)\right]-\lambda_na\int\limits_0^af_n(s)s^{\alpha_1-1}E_{(\alpha_1-\alpha_2,\alpha_1),\alpha_1}\left(-\mu s^{\alpha_1-\alpha_2},-\lambda_ns^{\alpha_1}\right)ds\},
$$
where,
\begin{equation}\label{con1}
	\Delta_n\equiv E_{(\alpha_1-\alpha_2,\alpha_1),1}\left(-\mu a^{\alpha_1-\alpha_2},-\lambda_na^{\alpha_1}\right)-1+\lambda_naE_{(\alpha_1-\alpha_2,\alpha_1),2}\left(-\mu a^{\alpha_1-\alpha_2},-\lambda_na^{\alpha_1}\right)\neq 0,
\end{equation}
Regarding the condition \eqref{con1}, we note that only for counted number of $n$ it can be violated. Moreover, one can show that $\lim\limits_{n\to\infty}|\Delta_n|>0$ and $\frac{1}{|\Delta_n|}<c$ ($c$ is a positive constant). 

If $1<\beta\le 4/3$, then $E_{\beta,2}(-z)$ for $z>0$ has no zeros \cite{R12}. Hence,  
\[
\Delta_n^*\equiv E_{\beta,2}\left(-\lambda_n(b-a)^\beta\right)>0    
\]
for all natural $n$ and one can show that $\lim\limits_{n\rightarrow\infty}\Delta_n\neq 0$, $\frac{1}{|\Delta_n^*|}<c$ ($c$ is a positive constant). Therefore, we find 
$$
G_n=1/[(b-a)\Delta_n^*]\{\psi_n-\int\limits_a^bf_n(s)(s-a)^{\beta-1}E_{\beta,\beta}\left(-\lambda_n(s-a)^\beta\right)ds
$$
$$
-E_{\beta,1}\left(-\lambda_n(b-a)^\beta\right)[\phi_n(1-\mu a^{\alpha_1-\alpha_2}E_{(\alpha_1-\alpha_2,\alpha_1),\alpha_1-\alpha_2+1}\left(-\mu a^{\alpha_1-\alpha_2},-\lambda_na^{\alpha_1}\right))
$$
$$
+\int\limits_0^af_n(s)s^{\alpha_1-1}E_{(\alpha_1-\alpha_2,\alpha_1),\alpha_1}\left(-\mu s^{\alpha_1-\alpha_2},-\lambda_ns^{\alpha_1}\right)ds]-B_nE_{(\alpha_1-\alpha_2,\alpha_1),2}\left(-\mu a^{\alpha_1-\alpha_2},-\lambda_na^{\alpha_1}\right)E_{\beta,1}\left(-\lambda_n(b-a)^\beta\right)\}.
$$

Now, we will use Lemma 2 to find an estimate for $ E_{(\alpha_1-\alpha_2,\alpha_1),\rho}\left(-\mu t^{\alpha_1-\alpha_2},-\lambda_nt^{\alpha_1}\right)$, but we first need to note the following property of the function \eqref{ML_mv}  (see \cite{R13}) before using the estimate given in Lemma 2:
\[
E_{(\alpha,\beta),\,\rho}(x,y)=E_{(\beta,\alpha),\,\rho}(y,x).
\]
Namely,
\[
E_{(\alpha_1-\alpha_2,\alpha_1),\rho}\left(-\mu t^{\alpha_1-\alpha_2},-\lambda_nt^{\alpha_1}\right)=E_{(\alpha_1,\alpha_1-\alpha_2),\rho}\left(-\lambda_nt^{\alpha_1}, -\mu t^{\alpha_1-\alpha_2}\right)
\]
and hence, according to Lemma 2 (note that $\alpha_1\in (1,2),\,\alpha_1>\alpha_2>0$ and $\rho>0,\,\mu>0,\,\lambda_n>0$), we have
\[
\left|E_{(\alpha_1-\alpha_2,\alpha_1),\rho}\left(-\mu t^{\alpha_1-\alpha_2},-\lambda_nt^{\alpha_1}\right)\right|\le \frac{c}{\lambda_nt^{\alpha_1}}.
\]

Considering, the above estimate, we can get the following estimates for $B_n, C_n$ and $D_n$:
\begin{eqnarray}\label{estcoef}
	&|B_n|\le \frac{1}{|\Delta_n|}\left[\frac{\mu c}{\lambda_n a^{\alpha_2}}|\phi_n|+c\lambda_n|\phi_n|+c |f_n(a)|+\lambda_na\int\limits_0^a |f_n(s)|s^{\alpha_1-1}ds\right], \nonumber\\
	&|C_n|\le \frac{c}{\lambda_na^{\alpha_1}}|\phi_n|+\frac{c}{\lambda_na^{\alpha_1}}|B_n|+\int\limits_0^a|f_n(s)|s^{\alpha_1-1}ds,\\
	&|G_n|\le \frac{1}{(b-a)|\Delta_n^*|}\left[|\psi_n|+\int\limits_a^b |f_n(s)|(s-a)^{\beta-1}ds+\frac{c}{\lambda_n(b-a)^\beta}\left[\frac{c}{\lambda_na^{\alpha_1}}|\phi_n|+\int\limits_0^a|f_n(s)|s^{\alpha_1-1}ds\right]+|B_n|\frac{c}{\lambda_n^2a^{\alpha_1}(b-a)^\beta}\right].\nonumber
\end{eqnarray}
Based on \eqref{Tn1_sol} and \eqref{Tn2_sol}, considering estimates of functions \eqref{ML_2}, \eqref{ML_mv}, one can obtain
\begin{eqnarray}\label{estTna}
	&|T_n(t)|\le c|\phi_n|+\frac{c}{\lambda_nt^{\alpha_1}}|B_n|+\int\limits_0^t|f_n(s)|s^{\alpha_1-1}ds,\,\,\,0< t\le a,\\
	&|T_n(t)|\le \frac{c}{\lambda_n(t-a)^\beta}|C_n|+\frac{c}{\lambda_n(t-a)^{\beta-1}}|G_n|+\int\limits_0^t |f_n(s)|(s-a)^{\beta-1}ds,\,\,\,a< t\le b.\label{estTnb}
\end{eqnarray}
Here we denote all positive constants via $c$, which have no principal influence. Now considering \eqref{estcoef}, we can rewrite \eqref{estTna} and \eqref{estTnb} as follows:
\begin{equation}\label{estTna2}
	|T_n(t)|\le c|\phi_n|+\frac{c}{t^{\alpha_1}|\Delta_n|}|\phi_n|+\frac{c}{\lambda_n^2t^{\alpha_1}|\Delta_n|}|\phi_n|+\frac{c}{\lambda_n t^{\alpha_1}|\Delta_n|}|f_n(a)|+\frac{c}{t^{\alpha_1}|\Delta_n|}\int\limits_0^a|f_n(s)|s^{\alpha_1-1}ds+\int\limits_0^t|f_n(s)|s^{\alpha_1-1}ds,\,\,0< t\le a,    
\end{equation}
\begin{eqnarray}\label{estTnb2}
	&|T_n(t)|\le \dfrac{c}{\lambda_n(t-a)^\beta|\Delta_n|}|\phi_n|+\dfrac{c}{\lambda_n(t-a)^{\beta-1}|\Delta_n^*|}|\psi_n|+\dfrac{c}{\lambda_n^2(t-a)^\beta}|\phi_n|+\dfrac{c}{\lambda_n^2(t-a)^{\beta-1}|\Delta_n^*|}|\phi_n|+\dfrac{c}{\lambda_n^3(t-a)^\beta|\Delta_n|}|\phi_n|\nonumber\\
	&+\dfrac{c}{\lambda_n^3(t-a)^{\beta-1}|\Delta_n^*|}|\phi_n|+\dfrac{c}{\lambda_n^4(t-a)^{\beta-1}|\Delta_n^*|}|\phi_n|+\dfrac{c}{\lambda_n^2(t-a)^\beta|\Delta_n|}|f_n(a)|+\dfrac{c}{\lambda_n^3(t-a)^{\beta-1}|\Delta_n^*|}|f_n(a)|\nonumber\\
	&+\dfrac{c}{\lambda_n(t-a)^\beta|\Delta_n|}\int\limits_0^a|f_n(s)|(s-a)^{\beta-1}ds+\dfrac{c}{\lambda_n^2(t-a)^\beta|\Delta_n|}\int\limits_0^a|f_n(s)|s^{\alpha_1-1}ds+\dfrac{c}{\lambda_n^2(t-a)^\beta}\int\limits_0^a|f_n(s)|s^{\alpha_1-1}ds\nonumber\\
	&+\dfrac{c}{\lambda_n(t-a)^{\beta-1}|\Delta_n^*|}\int\limits_a^b|f_n(s)|(s-a)^{\beta-1}ds+\dfrac{c}{\lambda_n^2(t-a)^{\beta-1}|\Delta_n^*|}\int\limits_0^a|f_n(s)|s^{\alpha_1-1}ds+\int\limits_0^t|f_n(s)|(s-a)^{\beta-1}ds,\,\,a<t\le b.
\end{eqnarray}

The convergence of infinite series $\sum\limits_{n=1}^\infty T_n(t)\sin(n\pi x)$ is based on the following estimates using \eqref{estTna2} and \eqref{estTnb2}:
\begin{eqnarray}\label{estTna2s}
	&\sum\limits_{n=1}^\infty|T_n(t)|\le c\left(1+\dfrac{1}{t^{\alpha_1}}\right)\sum\limits_{n=1}^\infty|\phi_n|+\dfrac{c}{t^{\alpha_1}}\sum\limits_{n=1}^\infty\dfrac{|\phi_n|}{\lambda_n^2}+\dfrac{c}{t^{\alpha_1}}\sum\limits_{n=1}^\infty\dfrac{|f_n(a)|}{\lambda_n}\nonumber\\
	&+\dfrac{c}{t^{\alpha_1}}\int\limits_0^a\sum\limits_{n=1}^\infty|f_n(s)|s^{\alpha_1-1}ds+\int\limits_0^t\sum\limits_{n=1}^\infty|f_n(s)|s^{\alpha_1-1}ds,\,\,0< t\le a,    
\end{eqnarray}
\begin{eqnarray}\label{estTnb2s}
	&\sum\limits_{n=1}^\infty|T_n(t)|\le \dfrac{c}{(t-a)^\beta}\sum\limits_{n=1}^\infty\dfrac{|\phi_n|}{\lambda_n}+\dfrac{c}{(t-a)^\beta}\sum\limits_{n=1}^\infty\dfrac{|\psi_n|}{\lambda_n}+\dfrac{c}{(t-a)^\beta}\sum\limits_{n=1}^\infty\dfrac{|\phi_n|}{\lambda_n^2}+\dfrac{c}{(t-a)^\beta}\sum\limits_{n=1}^\infty\dfrac{|\phi_n|}{\lambda_n^3}+\dfrac{c}{(t-a)^{\beta-1}}\sum\limits_{n=1}^\infty\dfrac{|\phi_n|}{\lambda_n^4}\nonumber\\
	&+\dfrac{c}{(t-a)^\beta}\sum\limits_{n=1}^\infty\dfrac{|f_n(a)|}{\lambda_n^2}+\dfrac{c}{(t-a)^{\beta-1}}\sum\limits_{n=1}^\infty\dfrac{|f_n(a)|}{\lambda_n^3}+\dfrac{c}{(t-a)^\beta}\int\limits_0^a\sum\limits_{n=1}^\infty\dfrac{|f_n(s)|}{\lambda_n}(s-a)^{\beta-1}ds+\dfrac{c}{(t-a)^\beta}\int\limits_0^a\sum\limits_{n=1}^\infty\dfrac{|f_n(s)|}{\lambda_n^2}s^{\alpha_1-1}ds\nonumber\\
	&+\dfrac{c}{(t-a)^{\beta-1}}\int\limits_a^b\sum\limits_{n=1}^\infty\dfrac{|f_n(s)|}{\lambda_n}(s-a)^{\beta-1}ds+\int\limits_0^t\sum\limits_{n=1}^\infty|f_n(s)|(s-a)^{\beta-1}ds,\,\,a<t\le b.
\end{eqnarray}
The convergence of all infinite series appearing in (\ref{estTna2s}) and (\ref{estTnb2s}) are stated in the following lemma:

{\bf Lemma 3.}
\begin{enumerate} [1)\,]
	\item	If $\phi(x)\in C[0,1]$ such that $\phi (0) = \phi (1) = 0$ and $\phi '\left( x \right) \in {L_2}(0,1),$ then \[\sum\limits_{n = 1}^\infty  {\left| {{\phi _n}} \right|}  \le \sum\limits_{n = 1}^\infty  {\frac{1}{{{{(n\pi )}^2}}} + \left\| {\phi '(x)} \right\|_2^2} ;\]
	
	\item	If $\phi(x)\in C^2[0,1]$ such that $\phi (0) = \phi (1) = 0,\phi ''(0) = \phi ''(1) = 0$ and $\phi '''\left( x \right) \in {L_2}(0,1),$ then
	\[\sum\limits_{n = 1}^\infty  {{{\left( {n\pi } \right)}^2}\left| {{\phi _n}} \right|}  \le \sum\limits_{n = 1}^\infty  {\frac{1}{{{{\left( {n\pi } \right)}^2}}} + \left\| {\phi '''} \right\|_2^2} ;\]
	
	\item	If $\psi(x)\in C[0,1]$ such that $\psi (0) = \psi (1) = 0$ and $\psi '\left( x \right) \in {L_2}(0,1),$ then \[\sum\limits_{n = 1}^\infty  {\left| {{\psi _n}} \right|}  \le \sum\limits_{n = 1}^\infty  {\frac{1}{{{{(n\pi )}^2}}} + \left\| {\psi '(x)} \right\|_2^2} ;\]

	\item	If $f(t,x)\in C[0,1]$ for all $t\in[0,T]$ such that $f\left( {t,0} \right) = 0,f\left( {t,1} \right) = 0$ and $\dfrac{{\partial f\left( {t,x} \right)}}{{\partial x}} \in {L_2}\left( {0,1} \right),$ then   \[\sum\limits_{n = 1}^\infty  {\left| {{f_n}\left( t \right)} \right|}  \le \sum\limits_{n = 1}^\infty  {\frac{1}{{{{\left( {n\pi } \right)}^2}}} + \left\| {\frac{{\partial f\left( {t,x} \right)}}{{\partial x}}} \right\|_2^2} ;\]
	
	\item	If $f(t,x)\in C^2[0,1]$ for all $t\in[0,T]$ such that $f\left( {t,0} \right) = 0,f\left( {t,1} \right) = 0,\,\dfrac{{{\partial ^2}f\left( {t,0} \right)}}{{\partial {x^2}}} = 0,\,\dfrac{{{\partial ^2}f\left( {t,1} \right)}}{{\partial {x^2}}} = 0,$ and $\dfrac{{{\partial ^3}f\left( {t,x} \right)}}{{\partial {x^3}}} \in {L_2}\left( {0,1} \right),$ then \[\sum\limits_{n = 1}^\infty  {{{\left( {n\pi } \right)}^{2}}\left| {{f_n}\left( t \right)} \right|}  \le \sum\limits_{n = 1}^\infty  {\frac{1}{{{{\left( {n\pi } \right)}^{2}}}} + \left\| {\frac{{{\partial ^3}f\left( {t,x} \right)}}{{\partial {x^3}}}} \right\|_2^2} ;\]
\end{enumerate}

The proof of this lemma can be done by using integration by parts, Bessel's inequality, and Parseval's identity.

For the uniform convergence of the series corresponding to the function $u_{xx}(t,x)$ we need to estimate $\sum\limits_{n=1}^\infty|\lambda_nT_n(t)|$. If one considers \eqref{estTna2s} and \eqref{estTnb2s}, one can easily see that stronger regularity will require the following terms:
\[
\sum\limits_{n=1}^\infty \lambda_n|\phi_n|,\,\,\sum\limits_{n=1}^\infty \lambda_n|f_n(s)|,\,\,\sum\limits_{n=1}^\infty |\psi_n|.
\] 

Regarding the uniqueness of the solution, we note that it follows from the completeness and basis property of the set of functions $\{\sin n\pi x\}_{n=1}^\infty$\cite{R14}.

Therefore, we can now formulate the main result in the following theorem:

{\bf Theorem 1.} {\it Let $1<\beta\le 4/3$ and \eqref{con1} be fulfilled. If conditions 2), 3), and 5) of Lemma 3 are valid, then problem \eqref{u_eqn}-\eqref{u_transm2} has a unique solution represented by \eqref{u_FS}.}

\vskip 2mm

\section*{\; Acknowledgment}
The authors acknowledge financial support from the National University of Science and Technology (NU), Oman. NU funds this work under research grant no. CFRG 23-02.

\end{document}